\documentclass{amsart}
\usepackage{amsfonts}

\usepackage{graphicx}

\input{tcilatex}
\input tcilatex

\begin{document}
\title[Zipf]{On divisibility properties of Zipf law and relatives}
\author{Thierry E. Huillet}
\address{Laboratoire de Physique Th\'{e}orique et Mod\'{e}lisation \\
CY Cergy Paris Universit\'{e}, CNRS UMR-8089,\\
2 avenue Adolphe-Chauvin, 95302 Cergy-Pontoise, FRANCE\\
E-mail: Thierry.Huillet@cyu.fr}
\maketitle

\begin{abstract}
We describe recursive chance mechanisms and processes involving Zipf's law
and relatives, all showing up mixing properties, both geometric and Poisson.
Relations to Markov processes are explored.\newline
\newline
\keywords{\textbf{Keywords:} Zipf laws, geometric and Poisson mixing,
self-decomposability, infinite divisibility, compound geometric and Poisson,
death process with immigration, Markov chain with disasters, heavy tails.\\
\\
AMS Classification: Primary: 60G, 60J, 60Exx, Secondary: 92D25.}
\end{abstract}

\section{Introduction}

The Zipf (or zeta) distribution is an ubiquitous discrete probability
distribution widely used in Sciences from various disciplines, including
linguistics, physics, genetics, sociology, fiability and reliability
engineering, among others. We herewith emphasize various divisibility
properties of this distribution and one of its companion version, focusing
on their Poisson and geometric mixing properties. As a result, they are
shown to appear as limiting laws of various Markov processes of interest.

The probability generating function (pgf) of Poisson-mixtures are in
one-to-one correspondence with the Laplace-Stieltjes transform (LSt) of the
mixing distribution, being absolutely monotone on the half-line with upper
bound their convergence radius. If in addition the mixing distribution is
self-decomposable (SD) [hence infinitely divisible], the associated pgf is
discrete self-decomposable (DSD) and unimodal [hence compound-Poisson or
discrete infinitely divisible], with a remarkable representation in terms of
the limit law of pure-death continuous-time branching processes with
immigration. A typical and central example is the pgf corresponding to a
variant of Zipf probability mass function (pmf), the mixing distribution of
which is derived here as a scale mixture of exponential random variables
(rv's). The classical Pareto-Zipf pmf shares the same properties with
modified mixing distribution. Such mixed Poisson distributions have
multiplicative explicit pmf's and they are heavy-tailed.

If in addition a pgf is a mixture of a geometric random variable (rv), the
associated pmf is completely monotone, hence log-convex. This is the case
for both Zipfian distributions described above, the mixing distributions
being related to log-Gamma distributions$.$ In that case, such distributions
are also compound geometric, a sub-class of compound-Poisson distributions.
As geometric mixtures, both Zipf's distributions enjoy the following
additional properties: Loss of memorylessness, decreasing failure rates and
overdispersion.

These distributions naturally appear in a recurrent Markovian random walk
with total resetting, respectively as its height of a typical excursion and
its stationary distribution. For this random walk (appearing in reliability
theory), depending on its positive or null recurrence, we further analyze
the number of returns to the origin (failures), the time to $k$-th return
and the contact probability at $0$ where Zipf distributions play a central
role.

A Zipf distribution also appears in a Bernoulli trial sequence (with
asymptotically harmonic) Zipfian unequal probabilities as the first return
time to $0$.

As a consequence of the Euler product formula for zeta function over the
subset of prime numbers, some additional divisibility properties of Zipf
distribution are finally mentioned.

\section{The Log-gamma model}

In this Section, we exhibit the strong divisibility properties of two
versions of Zipf's law.

\subsection{A variant of Zipf distribution, as a geometric mixture}

Let $\Gamma $ be Gamma$\left( \alpha ,\beta \right) -$distributed (with
density $f_{\Gamma }\left( x\right) =\frac{\alpha ^{\beta }}{\Gamma \left(
\beta \right) }x^{\beta -1}e^{-\alpha x}$) and $U=e^{-\Gamma }$ a
log-Gamma-distributed rv (with image density $f_{U}\left( u\right) =\frac{%
\alpha ^{\beta }}{\Gamma \left( \beta \right) }\left( -\log u\right) ^{\beta
-1}u^{\alpha -1}$). We have 
\begin{eqnarray*}
\int_{0}^{1}\left( -\log u\right) ^{\beta -1}u^{\alpha -1}du &=&\Gamma
\left( \beta \right) \alpha ^{-\beta } \\
\int_{0}^{1}u^{k}\left( -\log u\right) ^{\beta -1}u^{\alpha -1}du &=&\Gamma
\left( \beta \right) \left( \alpha +k\right) ^{-\beta } \\
\left[ z^{k}\right] \int_{0}^{1}\frac{1-u}{1-uz}\left( -\log u\right)
^{\beta -1}u^{\alpha -1}du &=&\Gamma \left( \beta \right) \left[ \left(
\alpha +k\right) ^{-\beta }-\left( \alpha +k+1\right) ^{-\beta }\right]
\end{eqnarray*}
\begin{equation}
\mathbf{P}\left( D=k\right) =:\pi _{k}=\left( 1+\frac{k}{\alpha }\right)
^{-\beta }-\left( 1+\frac{k+1}{\alpha }\right) ^{-\beta }  \label{pmf}
\end{equation}

\begin{equation*}
\mathbf{P}\left( D\geq k\right) =\left( 1+\frac{k}{\alpha }\right) ^{-\beta }%
\text{, }k\in \Bbb{N}_{0}=\left\{ 0,1,2,\cdots \right\}
\end{equation*}
\begin{equation}
\mathbf{E}\left( z^{D}\right) =\frac{\alpha ^{\beta }}{\Gamma \left( \beta
\right) }\int_{0}^{1}\frac{1-u}{1-uz}\left( -\log u\right) ^{\beta
-1}u^{\alpha -1}du  \label{pgf}
\end{equation}
showing that the rv $D$ with pmf (\ref{pmf}) is a geometric mixture with
mixing density $f_{U}\left( u\right) .$ Note that $\mathbf{P}\left( D\geq
k\right) \sim \left( \frac{k}{\alpha }\right) ^{-\beta }$. As a geometric
mixture, the pmf $\left( \pi _{k}\right) _{k}$ enjoy the complete
monotonicity property: 
\begin{equation*}
\left( -1\right) ^{n}\Delta ^{n}\pi _{k}\geq 0\text{ for all }n,k\in \Bbb{N}%
_{0}=\left\{ 0,1,2,\cdots \right\} ,
\end{equation*}
where $\Delta $ is the forward shift operator. In particular, as a Hausdorff
moment sequence, such pmf's are log-convex (obeying $\pi _{k}^{2}\leq \pi
_{k+1}\pi _{k-1}$, $\in \Bbb{N}=\left\{ 1,2,\cdots \right\} $).

In the sequel, we shall limit ourselves to the case when the scale parameter
is $\alpha =1.$ For such rv's $D$, for which $\pi _{k}=\left( 1+k\right)
^{-\beta }-\left( 2+k\right) ^{-\beta }\sim \beta k^{-\left( \beta +1\right)
}$ as $k\rightarrow \infty $, with $\zeta _{\beta }\left( z\right)
=\sum_{k\geq 1}k^{-\beta }z^{k}$ the polylog series, the pgf reads: 
\begin{eqnarray*}
\phi _{D}\left( z\right) &:&=\mathbf{E}\left( z^{D}\right) =\sum_{k\geq
0}\left( k+1\right) ^{-\beta }z^{k}-\sum_{k\geq 0}\left( k+2\right) ^{-\beta
}z^{k} \\
&=&z^{-1}\zeta _{\beta }\left( z\right) -z^{-2}\sum_{k\geq 2}k^{-\beta
}z^{k}=z^{-1}\left( 1-\left( \frac{1-z}{z}\right) \zeta _{\beta }\left(
z\right) \right)
\end{eqnarray*}
So, 
\begin{equation}
\phi _{D}\left( z\right) =z^{-1}\left( 1-\left( 1-z\right) \zeta _{\beta
}^{0}\left( z\right) \right) \text{, }\zeta _{\beta }^{0}\left( z\right)
:=\sum_{k\geq 0}\left( k+1\right) ^{-\beta }z^{k}.  \label{D1}
\end{equation}
With $D_{+}=D+1,$ the shifted version of $D$, 
\begin{equation*}
\phi _{D_{+}}\left( z\right) :=\mathbf{E}\left( z^{D_{+}}\right) =1-\left(
1-z\right) \zeta _{\beta }^{0}\left( z\right) .
\end{equation*}
Both $\phi _{D}\left( z\right) $ and $\phi _{D_{+}}\left( z\right) $ are
proper pgf's if and only if $\beta >0$, both with convergence radius $1.$
They are heavy-tailed with tail-index $\beta .$ Both have infinite mean if $%
\beta \in \left( 0,1\right] $, finite mean ($=\zeta _{\beta }\left( 1\right)
-1,$ $\zeta _{\beta }\left( 1\right) >1$ respectively) and infinite variance
if $\beta \in \left( 1,2\right] $, finite variance if $\beta >2.$

Note $\zeta _{\beta +1}\left( z\right) =\int_{0}^{z}\frac{\zeta _{\beta
}\left( z^{\prime }\right) }{z^{\prime }}dz^{\prime }$ and $\zeta _{\beta
+1}^{0}\left( z\right) =\frac{1}{z}\int_{0}^{z}\zeta _{\beta }^{0}\left(
z^{\prime }\right) dz^{\prime }$ entailing $\zeta _{\beta +1}\left( z\right)
<\zeta _{\beta }\left( z\right) $ and $\zeta _{\beta +1}^{0}\left( z\right)
<\zeta _{\beta }^{0}\left( z\right) .$

\textbf{Remark.} If $D$ is to represent the lifetime of some device until
breakdown occurs, as a geometric mixture, it is the superposition of a
continuum of geometric components with random breakdown probability $1-U.$

- If $\beta <1:$ the mixing density $f_{U}\left( u\right) \propto \left(
-\log u\right) ^{\beta -1}$ goes to $\infty $ when $u\rightarrow 1.$ The
population is saturated with individuals with a very small breakdown
probability and the first breakdown is expected to only occur after a long
time. The resulting tail of $D$ is very heavy, with $\mathbf{P}\left(
D>k\right) \sim k^{-\beta }$, obeying $\mathbf{E}\left( D\right) =\infty $.

- If $\beta =1:$ the mixing density collapse to $f_{U}\left( u\right) $
uniform. The breakdown parameter $1-U$ is uniform and all breakdown
probabilities are equiprobable. The final distribution of $D$ is
(Yule-Simon): $\mathbf{P}\left( D=k\right) =1/\left[ \left( k+1\right)
\left( k+2\right) \right] $, with tail index $1$, still with $\mathbf{E}%
\left( D\right) =\infty $.

- If $\beta >1:$ When $u\rightarrow 0,$ $f_{U}\left( u\right) $ goes to $%
\infty :$ there is a majority of individuals with large breakdown
probability and the first breakdown is expected to occur after a very short
time. The probability that $D=1$ or $D=2$ is thus significantly large.
Concomitantly, the mixing density goes to $0$ when $u\rightarrow 1$ and so,
only a minority of individuals has a small breakdown probability. However,
their time to breakdown is so long that it compensates and overwhelms the
rapid time to breakdown of the majority (the Reverse Lottery Paradox). The
resulting tail of $D$ remains heavy, still with $\mathbf{P}\left( D>k\right)
\sim k^{-\beta }$ but now with $\mathbf{E}\left( D\right) <\infty $.

The random walk approach with total resetting developed below will
corroborate this heuristics in which the most fragile components disappear
first, only leaving the most robust ones alive.

\subsection{The companion Zipf distribution is also Poisson mixture}

We herewith show that the rv $D$ also is a Poisson mixture and we compute
its random intensity. Poisson mixtures appear in reliability theory: a
system may fail not from wear-and-tear, but because it is hit by external
shocks. The number of shocks hitting the device typically follows a Poisson
distribution. When the intensity of the shock environment is random, the
number of shocks is a Poisson mixture.

Recall that (when $\alpha =1$) with $f_{U}\left( u\right) =\frac{1}{\Gamma
\left( \beta \right) }\left( -\log u\right) ^{\beta -1}$%
\begin{equation*}
\mathbf{P}\left( D=k\right) =\int_{0}^{1}\left( 1-u\right) u^{k}f_{U}\left(
u\right) du.
\end{equation*}
Observing that, with $\theta =\left( 1-u\right) /u>0$: $\left( 1-u\right)
u^{k}=\theta \int_{0}^{\infty }\frac{e^{-x}x^{k}}{k!}e^{-\theta x}dx$,
applying Fubini's theorem 
\begin{equation}
\mathbf{P}\left( D=k\right) =\int_{0}^{\infty }dx\frac{e^{-x}x^{k}}{k!}%
\int_{0}^{\infty }\theta e^{-\theta x}f_{\Theta }\left( \theta \right)
d\theta  \label{D2}
\end{equation}
where $f_{\Theta }\left( \theta \right) =\left( 1+\theta \right)
^{-2}f_{U}\left( 1/\left( 1+\theta \right) \right) $ is the image
probability density of $\Theta =U^{-1}\left( 1-U\right) .$ This shows that $%
D $ is a Poisson mixture with mixing probability density 
\begin{equation}
f_{X}\left( x\right) =\int_{0}^{\infty }\theta e^{-\theta x}f_{\Theta
}\left( \theta \right) d\theta ,  \label{D3}
\end{equation}
the Laplace transform of $\theta f_{\Theta }\left( \theta \right) .$ Note 
\begin{equation*}
\int_{0}^{\infty }f_{X}\left( x\right) dx=\int_{0}^{\infty }d\theta \cdot
\theta f_{\Theta }\left( \theta \right) \int_{0}^{\infty }e^{-\theta
x}dx=\int_{0}^{\infty }f_{\Theta }\left( \theta \right) d\theta =1.
\end{equation*}

The mixing density $f_{X}\left( x\right) $ being a scale mixture of
exponential rv's, by Thorin theorem on Generalized Gamma Convolution rv's 
\cite{Bond}, $f_{X}\left( x\right) $ is self-decomposable (SD), unimodal
[and hence infinitely divisible], together with $D$ being discrete
self-decomposable (DSD), unimodal (here with mode at the origin) [hence
compound Poisson], see \cite{Steu}. Let $X$ be the rv with density $%
f_{X}\left( x\right) $. We have 
\begin{equation*}
L_{X}\left( s\right) :=\mathbf{E}\left( e^{-sX}\right) =\int_{0}^{\infty }%
\frac{\theta }{\theta +s}f_{\Theta }\left( \theta \right) d\theta
\end{equation*}
Because, as a LSt $L\left( s\right) $, $s\geq 0$ is completely monotone on $%
\left[ 0,\infty \right) $ (that is, obeys: $\left( -1\right)
^{n}L_{X}^{\left( n\right) }\left( s\right) \geq 0$ for all $n\geq 0$), 
\begin{equation}
\phi _{D}\left( z\right) =z^{-1}\left( 1-\left( 1-z\right) \zeta _{\beta
}^{0}\left( z\right) \right) =\int_{0}^{\infty }e^{-x\left( 1-z\right)
}f_{X}\left( x\right) dx=L_{X}\left( 1-z\right) ,  \label{D4}
\end{equation}
as a scalable pgf, is absolutely monotone on the full half-line (that is,
obeys: $\phi _{D}^{\left( n\right) }\left( s\right) \geq 0$ for all $n\geq 0$%
, all $s\in $ $\left( -\infty ,0\right] $) and not only on $s\in \left(
0,1\right] $, with the one-to-one correspondence, \cite{KS1}:

\begin{equation*}
\phi _{D}\left( z\right) =L_{X}\left( 1-z\right) \text{ and }L_{X}\left(
s\right) =\phi _{D}\left( 1-s\right) .
\end{equation*}

\subsection{The L\'{e}vy measure of the mixing rv $X$}

$X$ being SD, it is in particular infinitely divisible. In view of the
L\'{e}vy-Khintchine representation of the exponential LSt, taking its
log-Laplace exponent and averaging, 
\begin{eqnarray*}
\frac{\theta }{\theta +s} &=&\exp \left( -\int_{0}^{\infty }\left(
1-e^{-sx}\right) x^{-1}e^{-\theta x}dx\right) , \\
-\log \left( \frac{\theta }{\theta +s}\right) &=&\int_{0}^{\infty }\left(
1-e^{-sx}\right) x^{-1}e^{-\theta x}dx \\
L_{X}\left( s\right) &=&\mathbf{E}_{\Theta }\exp \left( -\int_{0}^{\infty
}\left( 1-e^{-sx}\right) x^{-1}e^{-\theta x}dx\right) \\
&=&\exp \left( -\int_{0}^{\infty }\left( 1-e^{-sx}\right) \pi \left(
x\right) dx\right) .
\end{eqnarray*}
where $\pi \left( x\right) =x^{-1}k\left( x\right) $. Let us characterize $%
k\left( x\right) .$ We have

\begin{equation*}
\Psi _{X}\left( s\right) :=-\log L_{X}\left( s\right) =\int_{0}^{\infty
}\left( 1-e^{-sx}\right) x^{-1}k\left( x\right) dx=-\log \int_{0}^{\infty }%
\frac{\theta }{\theta +s}f_{\Theta }\left( \theta \right) d\theta
\end{equation*}
so that 
\begin{equation}
\Psi _{X}^{\prime }\left( s\right) =\int_{0}^{\infty }e^{-sx}k\left(
x\right) dx=\frac{\int_{0}^{\infty }\frac{\theta }{\left( \theta +s\right)
^{2}}f_{\Theta }\left( \theta \right) d\theta }{\int_{0}^{\infty }\frac{%
\theta }{\theta +s}f_{\Theta }\left( \theta \right) d\theta }.  \label{D5}
\end{equation}

The L\'{e}vy-measure of $X$ is then $\pi \left( x\right) dx=x^{-1}\mathcal{L}%
^{-1}\left( \Psi ^{\prime }\left( s\right) \right) \left( x\right) dx$
(involving a Laplace inverse $\mathcal{L}^{-1}$); it integrates $1\wedge x$.
As $x\rightarrow 0^{+}$, $\pi \left( x\right) \sim C\cdot x^{-1}$
(infinitely many microscopic jumps over any finite time interval) while as $%
x\rightarrow \infty $, $\pi \left( x\right) \sim C\cdot x^{-2}\left( \log
x\right) ^{\beta }$ (heavy-tailed large jumps attenuated by a logarithmic
factor).

Observing that $\Psi _{X}^{\prime }\left( s\right) $ is completely monotone, 
$\Psi _{X}^{\prime }\left( 1-z\right) $ is absolutely monotone.

To summarize:

\begin{proposition}
The rv $D$ (defined for $\beta >0$) with pmf (\ref{pmf}) is a geometric
mixture with log-Gamma mixing distribution $f_{U}\left( u\right) $.
Consequently, its pmf is completely monotone and log convex. It is also a
Poisson mixture with self-decomposable (so infinitely divisible) mixing
distribution $f_{X}\left( x\right) $ in (\ref{D2})(\ref{D3}), with Laplace
exponent characterized by $\Psi _{X}^{\prime }\left( s\right) $ in (\ref{D5}%
). As such, $D$ is DSD. The rv $D$ is thus compound geometric (and therefore
compound Poisson), with heavy-tailed law with tail-index $\beta $, unimodal
with mode at the origin.
\end{proposition}

\subsection{The distribution of $D$ as the limit law of a branching process
with immigration}

$\bullet $ From the DSD property of $D$:

\begin{equation*}
\phi _{D}\left( z\right) =\mathbf{E}\left( z^{D}\right) =L_{X}\left(
1-z\right) =e^{-\Psi _{X}\left( 1-z\right) }=e^{-r\int_{z}^{1}\frac{%
1-h\left( z^{\prime }\right) }{1-z^{\prime }}dz^{\prime }}
\end{equation*}
where $h\left( z\right) =\mathbf{E}\left( z^{I}\right) $ is the pgf of the
number $I$ of immigrants (with $h\left( 0\right) =0$) and $e^{-r\int_{z}^{1}%
\frac{1-h\left( z^{\prime }\right) }{1-z^{\prime }}dz^{\prime }}$ the
limiting pgf of a continuous-time pure birth branching process with
immigration events occurring at times of increase of a standard Poisson
process, see (\cite{vanH}, \cite{Fos}, \cite{Steutel}, p. 454). We thus
have: 
\begin{equation*}
\Psi _{X}^{\prime }\left( 1-z\right) =r\frac{1-h\left( z\right) }{1-z}
\end{equation*}
where $\frac{1-h\left( z\right) }{1-z}$ is the survival pgf of $I$, with$:$ $%
\left[ z^{n}\right] \frac{1-h\left( z\right) }{1-z}=\mathbf{P}\left(
I>n\right) .$ $\Psi _{X}^{\prime }\left( 1\right) =r$ fixes $r>0$. It can be
checked that $\Psi _{X}^{\prime }\left( 0\right) =\int_{0}^{\infty }\theta
^{-1}f_{\Theta }\left( \theta \right) d\theta =\infty $, so that $\Psi
_{X}^{\prime }\left( 0\right) =r\frac{1-h\left( z\right) }{1-z}\mid _{z=1}=r%
\mathbf{E}\left( I\right) =\infty .$ We have 
\begin{equation*}
\left[ z^{n}\right] \Psi _{X}^{\prime }\left( 1-z\right) =\int_{0}^{\infty
}e^{-\left( 1-z\right) x}k\left( x\right) dx=\frac{1}{n!}\int_{0}^{\infty
}x^{n}e^{-x}k\left( x\right) dx,
\end{equation*}
so that, when $n$ is large, 
\begin{equation*}
\mathbf{P}\left( I=n\right) \sim C\cdot n^{-2}\left( \log n\right) ^{\beta },
\end{equation*}
a heavy-tailed behavior mimicking the behavior for large $x$ of the L\'{e}vy
density of $X$: a large number of immigrants are expected to arrive in batch
at each immigration event.

$\bullet $ Because $D$ is a geometric mixture, it is also compound geometric
(and thus also compound Poisson), (see Th. 1.5.13 of \cite{VH} page $21$),
so with 
\begin{equation}
\phi _{D}\left( z\right) =\frac{\phi _{D}\left( 0\right) }{1-\left( 1-\phi
_{D}\left( 0\right) \right) \phi _{\Delta }\left( z\right) }  \label{D6}
\end{equation}
for some pgf 
\begin{equation}
\phi _{\Delta }\left( z\right) =\frac{1}{1-\phi _{D}\left( 0\right) }\left(
1-\frac{\phi _{D}\left( 0\right) }{\phi _{D}\left( z\right) }\right) \text{, 
}\phi _{\Delta }\left( 0\right) =0.  \label{D7}
\end{equation}
Then, as $n\rightarrow \infty $, $D_{n}\overset{d}{\rightarrow }D$ where $%
\left( D_{n}\right) _{n}$ is the following discrete-time homogeneous random
walk with catastrophes 
\begin{eqnarray*}
D_{n+1} &=&D_{n}+\Delta _{n+1}\text{ w.p. }1-\phi _{D}\left( 0\right) \\
D_{n+1} &=&0\text{ w.p. }\phi _{D}\left( 0\right) ,
\end{eqnarray*}
for some i.i.d. sequence $\left( \Delta _{n}\right) _{n}$ with common law
the one of $\Delta $.

\subsection{Tail-biased Zipf model}

Introduce a tail-biased (defined only if $\beta >1$) version $K$ of $D$
(Zipf-Mandelbrot): 
\begin{equation}
\pi _{k}^{*}:=\mathbf{P}\left( K=k\right) =\frac{\mathbf{P}\left( D\geq
k\right) }{\zeta _{\beta ,\alpha }}=\frac{\left( 1+\frac{k}{\alpha }\right)
^{-\beta }}{\zeta _{\beta ,\alpha }},k\in \Bbb{N}_{0}  \label{Z1}
\end{equation}
where $\sum_{k\geq 1}\left( k/\alpha \right) ^{-\beta }:=\zeta _{\beta
,\alpha }<\infty $ where $\zeta _{\beta ,\alpha }$ is the Hurwitz zeta
number. $K$ remains a geometric mixture in view of 
\begin{equation*}
\sum_{k\geq 0}z^{k}\mathbf{P}\left( D\geq k\right) =\frac{\alpha ^{\beta }}{%
\Gamma \left( \beta \right) \left( 1-z\right) }\int_{0}^{1}\frac{1}{1-uz}%
\left( -\log u\right) ^{\beta -1}u^{\alpha -1}du.
\end{equation*}
The complete monotonicity of a pmf indeed entails the complete monotonicity
of its complementary cumulative pmf (also known as the survival or tail
function).

When $\alpha =1:$ $\mathbf{P}\left( K=k\right) =\left( 1+k\right) ^{-\beta
}/\zeta _{\beta }$, $k\in \Bbb{N}_{0}$ where $\zeta _{\beta }=\sum_{k\geq
1}k^{-\beta }$ is the zeta number. Note that, with $\zeta _{\beta }\left(
z\right) =\sum_{k\geq 1}k^{-\beta }z^{k}$ the Dirichlet (or polylogarithm)
series 
\begin{equation}
\phi _{K}\left( z\right) :=\mathbf{E}\left( z^{K}\right) =z^{-1}\frac{\zeta
_{\beta }\left( z\right) }{\zeta _{\beta }\left( 1\right) }=\frac{\zeta
_{\beta }^{0}\left( z\right) }{\zeta _{\beta }\left( 1\right) }=:z^{-1}\phi
_{K_{+}}\left( z\right) .  \label{Z2}
\end{equation}
$\phi _{K}\left( z\right) $ is absolutely monotone on $\left( -\infty
,z_{c}=0\right] $ in view of $\zeta _{\beta }\left( z\right) \sim -\log
^{\beta }\left( -z\right) /\Gamma \left( \beta +1\right) $ as $z\rightarrow
-\infty .$

With $K_{+}=K+1:$ $\pi _{k}^{*,+}:=\mathbf{P}\left( K_{+}=k\right)
=k^{-\beta }/\zeta _{\beta }$, $k\in \Bbb{N}$ and $\phi _{K_{+}}\left(
z\right) :=\mathbf{E}\left( z^{K_{+}}\right) =\frac{\zeta _{\beta }\left(
z\right) }{\zeta _{\beta }\left( 1\right) }$ is the pgf of the shifted
version $K_{+}$ of $K.$

Both $K$ and $K_{+}$ have infinite mean if $\beta \in \left( 1,2\right] $,
finite mean ($=\zeta _{\beta -1}/\zeta _{\beta }-1,$ $\zeta _{\beta
-1}/\zeta _{\beta }>1$ respectively) and infinite variance if $\beta \in
\left( 2,3\right] $, finite variance if $\beta >3.$

The pgf $\phi _{K}\left( z\right) $ is scalable being absolutely monotone on 
$\left( -\infty ,1\right] $.

- With $f_{U}^{*}\left( u\right) =f_{U}\left( u\right) /\left( 1-u\right) $,
we have 
\begin{equation*}
\mathbf{P}\left( K=k\right) =\frac{1}{\zeta _{\beta }}\int_{0}^{1}u^{k}f_{U}%
\left( u\right) du=\frac{1}{\zeta _{\beta }}\int_{0}^{1}\left( 1-u\right)
u^{k}f_{U}^{*}\left( u\right) du.
\end{equation*}
So $K$ itself is a geometric mixture with modified mixing density $\frac{1}{%
\zeta _{\beta }}f_{U}^{*}\left( u\right) $.

- Observing that, with $\theta =\left( 1-u\right) /u>0$: $\left( 1-u\right)
u^{k}=\int_{0}^{\infty }\frac{e^{-x}x^{k}}{k!}e^{-\theta x}dx$, applying
Fubini's theorem 
\begin{equation*}
\mathbf{P}\left( K=k\right) =\frac{1}{\zeta _{\beta }}\int_{0}^{\infty }dx%
\frac{e^{-x}x^{k}}{k!}\int_{0}^{1}e^{-\theta x}f_{\Theta _{*}}\left( \theta
\right) d\theta
\end{equation*}
where 
\begin{equation}
f_{\Theta _{*}}\left( \theta \right) =\left[ \theta \left( 1+\theta \right)
\right] ^{-1}f_{U}\left( 1/\left( 1+\theta \right) \right)  \label{Z3}
\end{equation}
\begin{equation*}
=\frac{1}{\Gamma \left( \beta \right) }\left[ \theta \left( 1+\theta \right)
\right] ^{-1}\log ^{\beta -1}\left( 1+\theta \right)
\end{equation*}
is the image density of $\Theta =U^{-1}\left( 1-U\right) $ when $U\sim
f_{U}^{*}\left( u\right) .$ This shows that $K$ is a Poisson mixture with
mixing probability density 
\begin{equation}
f_{X_{*}}\left( x\right) =\frac{1}{\zeta _{\beta }}\int_{0}^{\infty }\theta
e^{-\theta x}f_{\Theta _{*}}\left( \theta \right) d\theta ,  \label{Z4}
\end{equation}
involving the Laplace transform of $\theta f_{\Theta _{*}}\left( \theta
\right) =\frac{1}{\Gamma \left( \beta \right) }\left( 1+\theta \right)
^{-1}\log ^{\beta -1}\left( 1+\theta \right) .$ Note 
\begin{equation*}
\int_{0}^{\infty }f_{X_{*}}\left( x\right) dx=\frac{1}{\zeta _{\beta }}%
\int_{0}^{\infty }f_{\Theta _{*}}\left( \theta \right) d\theta =\frac{1}{%
\Gamma \left( \beta \right) \zeta _{\beta }}\int_{0}^{\infty }\frac{\log
^{\beta -1}\left( 1+\theta \right) }{\theta \left( 1+\theta \right) }d\theta
=1,
\end{equation*}
in agreement (upon the change of variable $1+\theta =e^{x}$) with the
integral representation $\zeta _{\beta }=\frac{1}{\Gamma \left( \beta
\right) }\int_{0}^{\infty }\frac{x^{\beta -1}}{e^{x}-1}dx,$ (see \cite{Hav},
p. $60$)$.$

As for $D$, the distribution of $K$ appears as the limit law of a branching
process with immigration, obtained while substituting $\Psi _{X_{*}}\left(
s\right) $ to $\Psi _{X}\left( s\right) .$ We skip details.

To summarize:

\begin{proposition}
The rv $K$ (defined for $\beta >1$) with pmf (\ref{Z1}) is a geometric
mixture with modified log-Gamma mixing distribution $f_{U}^{*}\left(
u\right) /\zeta _{\beta }$. Consequently, its pmf is completely monotone and
log convex. It is also a Poisson mixture with self-decomposable (so
infinitely divisible) mixing distribution $f_{X_{*}}\left( x\right) $ in (%
\ref{Z3})(\ref{Z4}) and with Laplace exponent $\Psi _{X_{*}}\left( s\right) $%
. As such, $K$ is DSD. The rv $K$ is thus compound geometric (and therefore
compound Poisson), with heavy-tailed law with tail-index $\beta -1$,
unimodal with mode at the origin.
\end{proposition}

\textbf{Remarks.} As geometric mixtures, both $D$ and $K$ Zipf's
distributions enjoy the following additional properties:

- Loss of Memorylessness: If a single geometric distribution is uniquely
memoryless (with constant failure rate), a mixture of distinct geometric
distributions explicitly breaks this property.

- Decreasing Failure Rate (DFR): A mixture of distinct geometric components
always exhibits a decreasing hazard rate (failure rate). Over time, the
''components'' with higher success/failure probabilities drop out early,
leaving the components with longer survival times. As a result, the longer a
system survives, the less likely it is to fail in the next interval (the
`richer gets richer' paradigm).

- Inflated Variance (overdispersion): As geometric mixtures, they are
overdispersed compared to a single geometric baseline with the same mean.
The variance of the mixture absorbs both the inherent randomness of the
geometric processes and the variance of the mixing distribution itself.

\section{A random walk with total disasters}

Let $\beta >0.$ Consider the new discrete time-homogeneous Markov chain
(MCD) with disasters $K:=\left( K_{n};n\geq 0\right) $ with state-space $%
\Bbb{N}_{0}$ and non-homogeneous spatial transition probabilities
characterized by:

Given $K_{n}=k\in \Bbb{N}_{0}$, the increment of $K_{n}$ is 
\begin{equation}
\begin{array}{l}
+1\text{ with probability}:\text{ }p_{k}=\left( 1+\left( k+1\right)
^{-1}\right) ^{-\beta }=\left( \frac{k+2}{k+1}\right) ^{-\beta } \\ 
-k\text{ with probability}:\text{ }q_{k}=1-\left( 1+\left( k+1\right)
^{-1}\right) ^{-\beta }.
\end{array}
\label{0d}
\end{equation}
Equivalently, $K_{n+1}=\left( K_{n}+1\right) \mathbf{1}\left(
U_{n+1}>q_{K_{n}}\right) $, where $\left( U_{n}\right) _{n}$ is an
independent and identically distributed (iid) sequence of uniform rv's. This
MCD models a population, the size of which grows linearly, until a
disastrous event takes it abruptly to extinction (state $0$), a state from
which the population starts growing again at fresh. In a different context:
at the beginning, when the fortune of some investor is small, the chance to
increase it is not significant, being often subject to ruin; if, by chance
the level of the fortune raised quite high, then the chance to increase it
becomes larger and larger (the rich gets richer paradigm). Over time, the
level of richness oscillates between low and high, highness here translating
into thick tails. In fiability theory, the lifetime of some device grows
linearly daily until its breakdown occurs: devices which survived breakdown
for long periods are more likely to survive even longer.

This MCD random walk is not time-reversible since detailed balance does not
hold. This is a slightly different MCD than the one developed in \cite{GH}
where there was some additional control on $p_{0}.$

For large $k$, the probabilities $p_{k}\sim 1-\beta /k$ and $q_{k}\sim \beta
/k$ are asymptotically harmonic: the disasters become rare when the
population size becomes large, a fact responsible for the heavy-tailed
nature of the limiting population size.

We could think of adding a base parameter $\alpha >-1$ (giving $p_{k}=\left(
\alpha +\left( k+1\right) ^{-1}\right) ^{-\beta }$) but we will not detail
this model as it would deal with the Hurwitz zeta function instead of the
zeta function that we have met before when $\alpha =1$.

For general $\left( p_{k},q_{k}\right) $ MCD's as in (\ref{0d}), the
following holds:

- If $C_{1}\equiv \sum_{k=0}^{\infty }q_{k}<\infty ,$ $K_{n}$ is transient.

- If $C_{1}=\infty ,$ $K_{n}$ is recurrent.

- If $C_{1}=\infty $ and $C_{2}\equiv \sum_{k\geq
0}\prod_{l=0}^{k-1}p_{l}<\infty ,$ $K_{n}$ is positive recurrent.

- If $C_{1}=\infty $ and $C_{2}\equiv \sum_{k\geq
0}\prod_{l=0}^{k-1}p_{l}=\infty ,$ $K_{n}$ is null recurrent.

In our case study, we have:

$C_{1}=\infty $ because $\beta >0$ and $C_{2}\equiv \sum_{k\geq
0}\prod_{l=0}^{k-1}p_{l}=\sum_{k\geq 1}k^{-\beta }=\zeta _{\beta }<\infty $
if and only if $\beta >1.$

So the specific MCD chain (\ref{0d}) with $\beta >0$ is recurrent, null
recurrent if $\beta \in \left( 0,1\right] $ and positive recurrent if $\beta
>1.$ When null recurrent, the chain has a non-finitely summable invariant
measure, whereas, when positive recurrent, the chain has an invariant
probability measure.

The MCD (\ref{0d}) started at $K_{0}=0$ is made of infinitely many iid
excursions between times to consecutive visits to $0$ (the extinction
times). Letting $\tau _{0,0}\geq 1$ be the time to consecutive visits to $0$%
, the height of each excursion is $H=\tau _{0,0}-1\geq 0.$ An excursion has
height $H=0$ whenever, starting from $0$, the chain remains in $0$ in the
next step, an event with probability $q_{0}=1-2^{-\beta }.$ For such trivial
excursions, $\tau _{0,0}=1.$ The rv $D$ defined above will be shown to be
the lifetime of the device (so the height $H$) in such a renewal process.
The stationary distribution of the rv $K_{\infty }$ will then be identified
with the Zipf law of $K$.

\subsection{\textbf{Return time to the origin of }$K_{n}$}

For the first return times to $0$, we have 
\begin{equation}
\mathbf{P}\left( \tau _{0,0}=k+1\right) =q_{k}\prod_{l=0}^{k-1}p_{l}\text{, }%
k\geq 0.  \label{13d}
\end{equation}
Equivalently, $\mathbf{P}\left( \tau _{0,0}>k\right)
=\prod_{l=0}^{k-1}p_{l}=\left( k+1\right) ^{-\beta },$ with $\mathbf{E}%
\left( \tau _{0,0}\right) =\sum_{k\geq 0}\mathbf{P}\left( \tau
_{0,0}>k\right) =\zeta _{\beta }<\infty $ if and only if $\beta >1.$ When $%
\beta \in \left( 0,1\right] $, $\mathbf{E}\left( \tau _{0,0}\right) =\infty
. $

The pgf $\phi _{0,0}\left( z\right) :=\mathbf{E}\left( z^{\tau
_{0,0}}\right) $ of $\tau _{0,0}\geq 1$ is 
\begin{eqnarray*}
\phi _{0,0}\left( z\right) &:&=\sum_{k\geq
0}q_{k}z^{k}\prod_{l=0}^{k-1}p_{l}=\sum_{k\geq 0}z^{k}\left( 1-\left(
1+\left( k+1\right) ^{-1}\right) ^{-\beta }\right) \prod_{l=1}^{k}\left(
1+l^{-1}\right) ^{-\beta } \\
&=&\sum_{k\geq 0}z^{k}\left( 1-\left( 1+\left( k+1\right) ^{-1}\right)
^{-\beta }\right) \left( k+1\right) ^{-\beta }
\end{eqnarray*}

\begin{eqnarray*}
&=&\sum_{k\geq 1}z^{k}\left( \left( k+1\right) ^{-\beta }-\left( k+2\right)
^{-\beta }\right) \\
&=&1-\frac{1}{z}\left( 1-z\right) \zeta _{\beta }\left( z\right) =1-\left(
1-z\right) \zeta _{\beta }^{0}\left( z\right) .
\end{eqnarray*}

It involves the Li$_{\beta }$-polylog function: $\zeta _{\beta }\left(
z\right) :=\sum_{k\geq 1}z^{k}k^{-\beta },$ with $\zeta _{\beta }^{0}\left(
z\right) =z^{-1}\zeta _{\beta }\left( z\right) $ its shifted version. Note $%
\mathbf{E}\left( z^{\tau _{0,0}}\right) =\phi _{D_{+}}\left( z\right) $,
showing that $\tau _{0,0}\overset{d}{=}D_{+}$ the shifted Zipf's variant of $%
D.$

- When $\beta \in \left( 0,1\right] $, $\mathbf{E}\left( \tau _{0,0}\right)
=\infty .$

- When $\beta >1$, the mean persistence time is finite, with 
\begin{equation*}
\phi _{0,0}^{\prime }\left( 1\right) =\mathbf{E}\left( \tau _{0,0}\right)
=C_{2}=\frac{1}{\pi _{0}}=\zeta _{\beta }.
\end{equation*}

The pgf of the height $H\geq 0$ of an excursion is 
\begin{equation}
\mathbf{E}\left( z^{H}\right) =:\phi _{H}\left( z\right) =z^{-1}\mathbf{E}%
\left( z^{\tau _{0,0}}\right) =z^{-1}\left( 1-\left( 1-z\right) \zeta
_{\beta }^{0}\left( z\right) \right) ,  \label{H}
\end{equation}
because $H=\tau _{0,0}-1.$ Note $\mathbf{E}\left( z^{H}\right) =\phi
_{D}\left( z\right) ,$ showing that $H\overset{d}{=}D$ as from the rv
defined in (\ref{pmf})$.$

The function $\phi _{H}\left( z\right) $ is the pgf of a discrete Pareto rv $%
H\geq 1$ with tail-index $\beta >0$: $\mathbf{P}\left( H>k\right) =\left(
k+1\right) ^{-\beta }$, $k\in \Bbb{N}_{0}$. We have: 
\begin{eqnarray*}
\mathbf{P}\left( H=k\right) &:&=\rho _{k}=\left[ z^{k}\right] \phi
_{H}\left( z\right) \underset{k\rightarrow \infty }{\sim }\beta k^{-\left(
\beta +1\right) }\text{ if }\beta \neq 1 \\
\mathbf{P}\left( H>k\right) &=&\frac{1}{1+k}\text{ if }\beta =1\text{ }
\end{eqnarray*}
and for all $k\geq 1$%
\begin{equation}
\rho _{k+1}/\rho _{k}<1,\text{ }\rho _{1}=1-2^{-\beta }  \label{f13b}
\end{equation}
($\rho _{k}$ is monotone decreasing). As a result, $\rho _{1}$ is the
maximal value of the $\rho _{k}$: the pmf $\rho _{k}$ has its mode at $k=1.$

When $\beta =1$, an easy computation shows that

\begin{equation*}
\phi _{H}\left( z\right) =1+\frac{1}{z}\left( 1-z\right) \log \left(
1-z\right) .
\end{equation*}

An estimate of $\phi _{H}\left( z\right) $\emph{\ } near $z=1$ is: 
\begin{equation*}
\begin{array}{l}
\phi _{H}\left( z\right) \underset{z\downarrow 1}{\sim }1-\Gamma \left(
1-\beta \right) \left( 1-z\right) ^{\beta }\text{ if }0<\beta <1\text{ } \\ 
\phi _{H}\left( z\right) \underset{z\downarrow 1}{\sim }1-\zeta _{\beta
}\cdot \left( 1-z\right) +O\left( 1-z\right) ^{\beta }\text{ if }\beta >1 \\ 
\phi _{H}\left( z\right) \underset{z\downarrow 1}{\sim }1+\left( 1-z\right)
\log \left( 1-z\right) \text{ if }\beta =1
\end{array}
\end{equation*}

\subsection{\textbf{Contact probability at }$0$.}

The pgf of $\tau _{0,0}$ appears in the evaluation of the contact
probability at time $n$. We have 
\begin{equation*}
\sum_{n\geq 0}\mathbf{P}_{0}\left( K_{n}=0\right) z^{n}=\frac{1}{1-\phi
_{0,0}\left( z\right) },
\end{equation*}
the Green function at $0,$ see \cite{GH}.

- If $0<\beta <1$, 
\begin{equation*}
1-\phi _{0,0}\left( z\right) \underset{z\downarrow 1}{\sim }\Gamma \left(
1-\beta \right) \left( 1-z\right) ^{\beta },
\end{equation*}
showing by singularity analysis that, in the null recurrent case (algebraic
decay of the contact probability) 
\begin{equation*}
\mathbf{P}_{0}\left( K_{n}=0\right) \underset{n\uparrow \infty }{\sim }\frac{%
1}{\Gamma \left( \beta \right) \Gamma \left( 1-\beta \right) }n^{-\left(
1-\beta \right) }.
\end{equation*}
When $\beta =1-\varepsilon $ ($\varepsilon >0$ small), the constant in front
of $n^{-\left( 1-\beta \right) }$ vanishes like $\varepsilon $.

- When $\beta =1$, 
\begin{equation*}
1-\phi _{0,0}\left( z\right) \underset{z\downarrow 1}{\sim }\left(
1-z\right) \log \left( 1-z\right) ,
\end{equation*}
with logarithmic singularity, showing by singularity analysis, that 
\begin{equation*}
\mathbf{P}_{0}\left( K_{n}=0\right) \underset{n\uparrow \infty }{\sim }\frac{%
1}{\log n}.
\end{equation*}

- When $\beta >1$ (positive recurrence case), using $\phi _{H}\left(
z\right) \underset{z\downarrow 1}{\sim }1-\zeta _{\beta }\left( 1\right)
\left( 1-z\right) +O\left( 1-z\right) ^{\beta },$ 
\begin{equation*}
1-\phi _{0,0}\left( z\right) \underset{z\downarrow 1}{\sim }\left(
1-z\right) \zeta _{\beta }
\end{equation*}
showing, as required that $\mathbf{P}_{0}\left( K_{n}=0\right) \underset{%
n\uparrow \infty }{\rightarrow }\pi _{0}=1/\zeta _{\beta }.$ When $\beta
=1+\varepsilon $ ($\varepsilon >0$ small), using $\zeta _{\beta }\underset{%
\beta \uparrow 1}{\sim }\varepsilon ^{-1}$: $\pi _{0}\sim \varepsilon $,
just like when $\beta <1.$

\subsection{\textbf{Invariant measure}}

The stationary distribution of the MCD, if it exists, should solve: 
\begin{equation*}
p_{0}\pi _{0}=\sum_{k\geq 1}\pi _{k}q_{k}\text{ and }\pi _{k}=\pi
_{0}\prod_{l=0}^{k-1}p_{l}\text{, }k\in \Bbb{N},
\end{equation*}
thereby showing up a multiplicative nature.

$\bullet $ In the positive recurrent case ($\beta >1$), with $\pi _{0}=1/%
\mathbf{E}\left( \tau _{0,0}\right) =1/\zeta _{\beta }$ (from Kac theorem),
we also have 
\begin{equation}
\pi _{k}=\frac{\prod_{l=0}^{k-1}p_{l}}{\mathbf{E}\left( \tau _{0,0}\right) }%
=\pi _{0}\prod_{l=1}^{k}\left( 1+l^{-1}\right) ^{-\beta }=\pi _{0}\left(
k+1\right) ^{-\beta }\text{, }k\in \Bbb{N}_{0},  \label{SD}
\end{equation}
so with 
\begin{equation*}
\mathbf{E}\left( z^{K_{\infty }}\right) =\zeta _{\beta }^{0}\left( z\right)
/\zeta _{\beta }\left( 1\right) .
\end{equation*}
This gives an explicit expression of the pgf of $K_{\infty }$ (in the
positive recurrent case) in terms of a shifted polylog function$.$ Note $%
\mathbf{E}\left( z^{K_{\infty }}\right) =\phi _{K}\left( z\right) $ so that $%
K_{\infty }\overset{d}{=}K$, the primary Zipf's rv.

$\bullet $ When $\beta \in \left( 0,1\right] $, $\pi _{k}=\pi _{0}\left(
k+1\right) ^{-\beta }$, $k\in \Bbb{N}_{0},$ for some arbitrary constant $\pi
_{0}$ and for which $\sum_{k\geq 0}\pi _{k}=\infty $ (non-summability of
such invariant measures).

\subsection{\textbf{Renewal aspects}}

In the positive recurrent case, we can also reproduce the arguments related
to renewal theory, for which 
\begin{equation*}
\Delta =\tau _{0,0}^{0}+\tau _{0,0}^{+},
\end{equation*}
where $\tau _{0,0}^{0}\overset{d}{\sim }$Geo$\left( p_{0}\right) $\ (viz $%
\mathbf{P}\left( \tau _{0,0}^{0}=k\right) =p_{0}q_{0}^{k}$, $k\geq 0$) is
the idle period (the number of consecutive trivial length-$1$ excursions)
and $\tau _{0,0}^{+}\geq 2$\ the length of a subsequent true excursion with
length $\geq 2$, for which: $\mathbf{P}\left( \tau _{0,0}^{+}=k+1\right)
=q_{k}\prod_{l=1}^{k-1}p_{l}$, $k\geq 1$, the lengths $\tau _{0,0}^{0}$, $%
\tau _{0,0}^{+}$ being mutually independent. Clearly, when $\beta >1$, with $%
C_{2}=\sum_{k\geq 0}\prod_{l=1}^{k-1}p_{l}=\zeta _{\beta }<\infty $%
\begin{equation*}
\mathbf{E}\left( \tau _{0,0}^{0}\right) =p_{0}/q_{0}\text{ and }\mathbf{E}%
\left( \tau _{0,0}^{+}\right) =:\mu ^{+}=1+C_{2}=1+\zeta _{\beta }>2.
\end{equation*}
gives the average contribution to $\Delta $ of the two components in the
special MCD case. The $\Delta $s now constitute the interval lengths of an
alternating renewal process, namely, the first time, starting from $0$, that
the walker moves back to $0$ from above.

\subsection{Number of returns to $0$}

Let $S_{n}=\sum_{m=1}^{n}\mathbf{1}\left( K_{m}=0\right) $, $n\geq 1$ and $%
T_{k}=\inf \left\{ n\geq 1:S_{n}=k\mid S_{0}=0\right\} $. With $\overline{%
\tau }_{k}:=\sum_{l=1}^{k}\tau _{0,0}^{\left( l\right) }$ the sum of the $k$
first return times to $0$, we have 
\begin{eqnarray*}
\mathbf{P}\left( S_{n}=k\right) &=&\mathbf{P}\left( \overline{\tau }_{k}\leq
n,\overline{\tau }_{k+1}>n\right) \\
&=&\mathbf{P}\left( \overline{\tau }_{k}\leq n\right) -\mathbf{P}\left( 
\overline{\tau }_{k+1}\leq n\right) \\
&=&\mathbf{P}\left( \overline{\tau }_{k+1}>n\right) -\mathbf{P}\left( 
\overline{\tau }_{k}>n\right) .
\end{eqnarray*}
\begin{equation*}
\mathbf{P}\left( S_{n}=k\right) =\left[ z^{n}\right] \frac{\phi _{0,0}\left(
z\right) ^{k}\left( 1-\phi _{0,0}\left( z\right) \right) }{1-z}\text{, }%
n\geq k.
\end{equation*}
The event $\left\{ T_{k}=n\right\} $ is realized whenever $n$ is the
end-point of the $k-$th return time to $0$. Hence 
\begin{equation*}
\mathbf{P}\left( T_{k}=n\right) =\mathbf{P}\left( \overline{\tau }%
_{k}=n\right) ,
\end{equation*}
with pgf 
\begin{equation*}
\mathbf{E}\left( z^{T_{k}}\right) =\phi _{0,0}\left( z\right) ^{k}.
\end{equation*}
As resulting from $\mathbf{P}\left( \tau _{0,0}>l\right) \sim l^{-\beta }$
as $l\rightarrow \infty $: 
\begin{equation*}
\mathbf{P}\left( \overline{\tau }_{k}>l\right) =\mathbf{P}\left(
T_{k}>l\right) \sim k\cdot l^{-\beta },
\end{equation*}
showing that $k^{-1/\beta }T_{k}$ is tail-equivalent to $\tau _{0,0}$.

\section{A related Bernoulli trial sequence with Zipfian unequal
probabilities}

Consider the independent Bernoulli sequence $\left( I_{\beta }\left(
l\right) \right) _{l}$ with $\mathbf{P}\left( I_{\beta }\left( l\right)
=0\right) =q_{l}.$ Assuming $\Sigma _{0}=0$, let $\Sigma _{n}=\sum_{l=1}^{n}%
\mathbf{1}\left( I_{\beta }\left( l\right) =0\right) $, $n\geq 1$, be the
number of occurrences of symbol $0$ of the Bernoulli string by time $n.$ The
sum process $\Sigma _{n}$ now is with independent and non-stationary
increments. We have

\begin{equation*}
\mathbf{E}\left( z^{\Sigma _{n}}\right) =\prod_{l=1}^{n}\left(
p_{l}+q_{l}z\right)
\end{equation*}
with $\mathbf{E}\left( \Sigma _{n}\right) =\sum_{l=1}^{n}q_{l}$ and $\sigma
^{2}\left( \Sigma _{n}\right) =\sum_{l=1}^{n}p_{l}q_{l}$. In our case study, 
$\mathbf{E}\left( \Sigma _{n}\right) \sim \beta \log n+O\left( 1\right) $
and $\sigma ^{2}\left( \Sigma _{n}\right) \sim \beta \log n$ with $\mathbf{E}%
\left( \Sigma _{n}\right) -\sigma ^{2}\left( \Sigma _{n}\right)
=\sum_{l=1}^{n}q_{l}^{2}\ll \mathbf{E}\left( \Sigma _{n}\right) $ and $%
\sum_{l=1}^{n}p_{l}q_{l}/\left( \log n\right) ^{2}<\infty $. $\Sigma _{n}$
is thus close, in total variation to a Poisson rv Poi$\left( \mathbf{E}%
\left( \Sigma _{n}\right) \right) $ so that $\left( \Sigma _{n}-\mathbf{E}%
\left( \Sigma _{n}\right) \right) /\sigma \left( \Sigma _{n}\right) $
converges in distribution to the reduced normal rv $\mathcal{N}\left(
0,1\right) $, together with $\Sigma _{n}/\left( \beta \log n\right)
\rightarrow 1$, a.s., as $n\rightarrow \infty .$ With $0\leq k\leq n$, we
have 
\begin{equation*}
\left[ z^{k}\right] \mathbf{E}\left( z^{\Sigma _{n}}\right) =\mathbf{P}%
\left( \Sigma _{n}=k\right) =\sum_{1\leq l_{1}<\cdots <l_{k}\leq
n}\prod_{j=1}^{k}q_{l_{j}}\prod_{l\neq \left\{ l_{1},\cdots ,l_{k}\right\}
}p_{l}.
\end{equation*}
Note the boundary probabilities: $\mathbf{P}\left( \Sigma _{n}=0\right)
=\prod_{l=1}^{n}p_{l}$, $\mathbf{P}\left( \Sigma _{n}=n\right)
=\prod_{l=1}^{n}q_{l}$ and 
\begin{equation*}
\mathbf{P}\left( \Sigma _{n}=1\right) =\sum_{k=1}^{n}q_{k}\prod_{l\neq
k}p_{l}=\sum_{k=1}^{n}\frac{q_{k}}{p_{k}}\prod_{l=1}^{n}p_{l}.
\end{equation*}
With $\mathcal{T}_{k}$ the time to $k$-th occurrence of $0$, the event $%
\Sigma _{n}=1$ is also the event $\mathcal{T}_{1}\leq n<\mathcal{T}_{2}$.

This probability sequence can be generated recursively according to 
\begin{equation*}
\mathbf{P}\left( \Sigma _{n+1}=k\right) =q_{n+1}\mathbf{P}\left( \Sigma
_{n}=k-1\right) +p_{n+1}\mathbf{P}\left( \Sigma _{n}=k\right) .
\end{equation*}
In contrast with (\ref{0d}), the transition probabilities of $\Sigma _{n}$
depend on time but they do not depend on the current state $k$.

Let $\mathcal{T}_{k}=\inf \left\{ n\geq 1:\Sigma _{n}=k\mid \Sigma
_{0}=0\right\} $ be time of the $k$-th occurrence of symbol $0$ (the inverse
process of $\Sigma _{n}$). With $n\geq k$, we have

\begin{equation*}
\mathbf{P}\left( \mathcal{T}_{k}=n\right) =q_{n}\mathbf{P}\left( \Sigma
_{n-1}=k-1\right)
\end{equation*}

We note that $\mathbf{P}\left( \mathcal{T}_{1}=n\right) =q_{n}\mathbf{P}%
\left( \Sigma _{n-1}=0\right) =q_{n}\prod_{l=1}^{n-1}p_{l}$, which coincides
with the time to first return to state $0$ of the MCD (\ref{0d})$.$ The
following recurrence can be checked to hold:

\begin{equation*}
\mathbf{P}\left( \mathcal{T}_{k+1}=n+1\right) =q_{n+1}\left( \mathbf{P}%
\left( \mathcal{T}_{k}=n\right) +\frac{p_{n}}{q_{n}}\mathbf{P}\left( 
\mathcal{T}_{k+1}=n\right) \right) ,
\end{equation*}
relating $\mathcal{T}_{k+1}$ to $\mathcal{T}_{k}.$

\section{Zeta law and primes, \protect\cite{Lin}}

Consider the positive zeta rv $K_{+}$ with values in $\Bbb{N}$. Its moment
function is 
\begin{equation*}
\mathbf{E}\left( K_{+}^{q}\right) =\frac{1}{\zeta _{\beta }}\sum_{k\geq 1}%
\frac{1}{k^{\beta -q}}=\frac{\zeta _{\beta -q}}{\zeta _{\beta }},\text{ for }%
q<\beta -1.
\end{equation*}
Thanks to the Euler product formula for zeta function over the subset of
primes $\mathcal{P}_{+}=\left\{ 2,3,5,\cdots \right\} $%
\begin{equation*}
\zeta _{\beta }=\prod_{p\in \mathcal{P}_{+}}\left( 1-p^{-\beta }\right)
^{-1},
\end{equation*}
we get 
\begin{equation*}
K_{+}\overset{d}{=}\prod_{p\in \mathcal{P}_{+}}p^{\text{Geo}\left( p^{-\beta
}\right) }\text{,}
\end{equation*}
where Geo$\left( p^{-\beta }\right) $, $p\in \mathcal{P}_{+}$ is a sequence
of independent shifted Geometric rv's with $\mathbf{P}\left( \text{Geo}%
\left( p^{-\beta }\right) =k\right) =p^{-\beta k}\left( 1-p^{-\beta }\right) 
$, $k\in \Bbb{N}_{0}$. The pgf of 
\begin{equation*}
\log K_{+}=\sum_{p\in \mathcal{P}_{+}}\text{Geo}\left( p^{-\beta }\right)
\log p
\end{equation*}
therefore is 
\begin{equation*}
\mathbf{E}\left( z^{\log K_{+}}\right) =\prod_{p\in \mathcal{P}_{+}}\frac{%
1-p^{-\beta }}{1-z^{\log p}p^{-\beta }},
\end{equation*}
an infinite product of shifted Geo$\left( p^{-\beta }\right) $ pgf's each
taking values on multiples of $\log p$, $p\in \mathcal{P}_{+},$ \cite{Lin}.
Note $\mathbf{P}\left( \log K_{+}=0\right) =\mathbf{P}\left( K_{+}=1\right) =%
\mathbf{P}\left( K=0\right) =1/\zeta _{\beta }$, consistently with the
analysis of $K$. Each such rv is DSD and so is the infinite product: the rv $%
\log \left( K_{+}\right) $ is DSD, an additional remarkable divisibility
feature of $K_{+}$.\newline

\textbf{Statements and Declarations.} The author declares that he has no
known competing financial interests or personal relationships that could
have appeared to influence the work reported in this paper.

\textbf{Conflicts of Interest / Competing Interests.}

The author has no conflicts of interest to declare that are relevant to the
content of this article.

\textbf{Acknowledgments.}

The author acknowledges partial support both from the Chair: ``\textit{%
Mod\'{e}lisation math\'{e}matique et biodiversit\'{e}'' }of\textit{\ }%
Veolia-Ecole Polytechnique-MNHN-Fondation X\textit{.}


\begin{thebibliography}{99}
\bibitem{Bond}  L. Bondesson, $\emph{Generalized}$ $\emph{Gamma}$\textrm{\ }$%
\emph{Convolutions}$\textrm{\ }$\emph{and}$\textrm{\ }$\emph{Related}$%
\textrm{\ }$\emph{Classes}$\textrm{\ }$\emph{of}$\textrm{\ }$\emph{%
Distributions}$\textrm{\ }$\emph{and}$\textrm{\ }$\emph{Densities.}$ Lecture
Notes in Statistics, vol. 76. Springer-Verlag, New York, 1992.

\bibitem{Bondesson}  L. Bondesson, A class of probability distributions that
is closed with respect to addition as well as multiplication of independent
random variables. \emph{J. Theor. Probab.} \textbf{28}, (2015), 1063-1081.

\bibitem{Fos}  E.B. Fosam, Characterization of discrete self-decomposable
distributions. \emph{Sankhya: The Indian Journal of Statistics}, \textbf{57}%
, Series A, Pt 3, (1995), 337-341.

\bibitem{GH}  B. Goncalves, T. Huillet. Scaling features of two Special
Markov chains involving total disasters. \emph{Journal of Statistical
Physics,} \textbf{178}(2), (2019), 499-531.

\bibitem{Hav}  J. Havil, The Gamma Function. Ch. 6 in Gamma: Exploring
Euler's Constant. Princeton, NJ: Princeton University Press, (2003), pp.
53-60.

\bibitem{KS1}  L. B. Klebanov, M. Sumbera, Discrete probability
distributions and statistical mechanics of small systems, (2026),
arXiv:2607.18968 [cond-mat.stat-mech].

\bibitem{Lin}  G.D. Lin, C-Y Hu, The Riemann zeta distribution. Bernoulli,
Volume 7, Number 5 (2001), 817-828.

\bibitem{Steutel}  F.W. Steutel, \& K. van Harn, \emph{Infinite Divisibility
of Probability Distributions on the Real Line,} Chapman \& Hall/CRC Pure and
Applied Mathematics. (2004).

\bibitem{Steu}  F.W. Steutel, \& K. van Harn, Discrete analogues of
self-decomposability and stability, \emph{Ann. Prob.,} \textbf{7} (1979),
893-899.

\bibitem{VH}  K. van Harn, \emph{Classifying infinitely divisible
distributions by functional equations}, Mathematical Centre Tracts, 103.
Mathematisch Centrum, Amsterdam, (1978), 194 pp.

\bibitem{vanH}  K. van Harn, F.W. Steutel, \& W. Vervaat, Self-decomposable
discrete distributions and branching processes, \emph{Z. Wahrsch. Verw.
Gebiete,} \textbf{61} (1982), 97-118.
\end{thebibliography}
\end{document}